 \documentclass[journal]{IEEEtran}

\ifCLASSINFOpdf
\else
\fi

\usepackage{svg}
\usepackage{xcolor}
\usepackage{amssymb}
\usepackage[mode=buildnew]{standalone}
\usepackage{enumerate}
\usepackage{amsmath}
\usepackage{mathtools}
\usepackage{mathrsfs}
\usepackage{tabularx}
\usepackage{comment}
\usepackage{graphicx}
\usepackage{psfrag}
\usepackage{epsfig}
\usepackage{multicol}
\usepackage{multirow}
\usepackage{graphics}
\usepackage{subcaption}
\usepackage{algorithm}
\usepackage{algcompatible}
\usepackage{booktabs}
\usepackage{float}
\usepackage{arydshln}
\usepackage{multirow}
\usepackage{multicol}
\usepackage{graphicx}
\usepackage[hidelinks]{hyperref}
\usepackage{cite}
\usepackage{lipsum}
\usepackage{optidef}

\usepackage{color}

\DeclareMathOperator{\sign}{sign}

\newtheorem{assumption}{Assumption}

\newtheorem{lemma}{Lemma}
\newtheorem{theorem}{Theorem}

 \newcommand\onenorm[1]{\left\lvert#1\right\rvert}

\usepackage{physics}

\begin{document}
\title{Extremum Seeking Control for Scalar Maps with Distributed Diffusion PDEs}    
	
	\author{Pedro Henrique Silva Coutinho, Tiago Roux Oliveira, and Miroslav Krsti\'{c}. 
		\thanks{P. H. S. Coutinho and T. R. Oliveira are with the Department of Electronics and Telecommunication Engineering,
			 State University of Rio de Janeiro (UERJ), Brazil
			(e-mails: phcoutinho@eng.uerj.br; tiagoroux@uerj.br).}
        \thanks{M. Krsti\'{c} is with the  Department of Mechanical and Aerospace Engineering, University of California - San Diego, USA (e-mail: krstic@ucsd.edu).}
		\thanks{Corresponding author: Pedro H. S. Coutinho. 
            }
		\thanks{This work was supported by the Brazilian agencies CNPq (Grant numbers: 407885/2023-4 and 309008/2022-0), CAPES and FAPERJ.}}

	\markboth{}%
	{Shell \MakeLowercase{\textit{et al.}}: Bare Demo of IEEEtran.cls for IEEE Journals}

	\maketitle

	\begin{abstract}
            This paper deals with the gradient extremum seeking control for static scalar maps with actuators governed by distributed diffusion partial differential equations (PDEs). To achieve the real-time optimization objective, we design a compensation controller for the distributed diffusion PDE via backstepping transformation in infinite dimensions. A further contribution of this paper is the appropriate motion planning design of the so-called probing (or perturbation) signal, which is more involved than in the non-distributed counterpart. Hence, with these two design ingredients, we provide an averaging-based methodology that can be implemented using the gradient and Hessian estimates. Local exponential stability for the closed-loop equilibrium of the average error dynamics is guaranteed through a Lyapunov-based analysis. By employing the averaging theory for infinite-dimensional systems, we prove that the trajectory converges to a small neighborhood surrounding the optimal point. The effectiveness of the proposed extremum seeking controller for distributed diffusion PDEs in cascade of nonlinear maps to be optimized is illustrated by means of numerical simulations.
	\end{abstract}


	\IEEEpeerreviewmaketitle

\section{Introduction}
\sloppy

Extremum-seeking control (ESC) is a real-time and model-free optimization approach
which is employed to determine and maintain the optimum (extremum) point of the output map
of a static or dynamic system with nonlinear and unknown dynamics \cite{scheinker2024100}.
For that purpose, the perturbation method \cite{Krsti2000StabilityOE} is often employed in the ESC.

ESC iteratively or adaptively adjusts the control-input variables in order to decrease or increase the value of a cost or objective function. This is done so that the system's output is gradually moved towards the optimal point, which can be a maximum or minimum of this objective function. 
ESC is particularly useful when it is not possible to obtain an accurate 
model of the system or when the system is subject to 
disturbances, uncertainties, or unpredictable variations~\cite{scheinker2024100}.

\paragraph*{Literature Review}
The ESC idea was initially proposed by Leblanc~\cite{leblanc1922} 
in the historical context of maximizing the power transfer to a tram car. 
Decades after, Krsti\'{c} and Wang in 
\cite{Krsti2000StabilityOE} have proposed a systematic analysis based on averaging and singular perturbation theories to prove stability for general nonlinear dynamic systems described by ordinary differential equations (ODEs). Since then, ESC have been improved in both practical and theoretical frameworks \cite{Ariyur2003RealTimeOB}.
To name a few examples, reference \cite{Manzie2009ExtremumSW} addresses the stochastic versions of ESC, the authors in 
\cite{Ghaffari2011MultivariableNE} propose the Newton-based ESC in order to obtain multiparameter generalization, while publication \cite{scheinker2017model} employs ESC for stabilization.
Moreover, ESC has been considered in several applications, 
such as brake system control, autonomous vehicles and mobile robots, 
bio-processes optimization, renewable energy systems, 
and refrigeration system optimization, see ~\cite{bastin2009extremum,tan2010extremum,nesic2012framework,krstic2014extremum,matveev2015extremum} and references therein. However, all this literature is limited to systems with finite-dimensional actuation dynamics.

It is known that actuation dynamics modeled by infinite-dimensional dynamics can be found in
several applications. For instance, actuator delays in which the delay propagation
is modeled by a first-order hyperbolic partial differential equation (PDE). 
Another example of an infinite-dimensional system is the diffusion process
that is found in biological, chemical, and temperature regulation~\cite{friedman2012pde,pearson2013fast,ng2012optimal}.

Within the context of ESC to infinite-dimensional actuation dynamics,
reference \cite{Oliveira2017ExtremumSF} addresses ESC 
with arbitrarily large delays whose effect is modeled by 
hyperbolic transport PDEs. This approach has paved the way for further generalizations to other classes of PDEs \cite{Oliveira_Krstic_book_2022}.
The authors in \cite{feiling2018gradient} handle the ESC with actuation dynamics governed by diffusion PDEs.
More recently, the Newton-based extremum seeking extension for
multi-input maps subject to actuation dynamics governed by distinct PDEs is performed in \cite{oliveira2020multivariable}, whereas publication \cite{silva2023extremum} considers the ESC for a class of wave PDEs with Kelvin-Voigt damping. 
However, the aforementioned results on ESC for PDEs are limited
to PDE boundary compensation. Those results can not be directly applied
when distributed actuation is concerned. \textcolor{black}{The sole ESC publication considering distributed delays is the paper \cite{Tsubakino_2023}.} Hence, the main goal of this paper is to expand the ESC framework for dealing with actuation dynamics governed
by distributed diffusion PDEs.

\paragraph*{Contributions}
Motivated by the aforementioned discussion, this paper deals with the
gradient ESC for static scalar maps with distributed actuation dynamics governed by 
diffusion PDEs. In contrast to related works~\cite{feiling2018gradient} and \cite{oliveira2020multivariable},
the boundary PDE compensation can not be applied in the distributed case. 
Then, inspired by the backstepping-like method proposed by~\cite{bekiaris2011compensating},
the first contribution of this paper is to develop a compensation controller for the distributed diffusion PDE
to guarantee the local exponential stability of the equilibrium of the closed-loop average error dynamics, based on
Lyapunov stability arguments.
Moreover, it is known that a key point for the successful operation of the ESC for PDE dynamics is the adequate design
of the so-called probing or perturbation signal. This motion planning problem is addressed by~\cite{feiling2018gradient,oliveira2020multivariable} by solving a 
trajectory generation with Dirichlet boundary condition~\cite{Krsti2008BoundaryCO}.
However, the perturbation signal designed in those works can not be directly employed when
distributed actuation is concerned. Therefore, the second contribution of this paper is to
solve a modified trajectory generation problem to design an adequate perturbation signal for
the ESC framework. With these two ingredients, we provide a novel ESC for distributed
diffusion PDE compensation which can be implemented using averaging-based estimates of the gradient and Hessian (first and second derivatives of the static map). 
Then, based on the exponential stability guarantees of the average error dynamics mentioned above, we
invoke the averaging theory for infinite-dimensional systems\cite{Hale1990AveragingII} to ensure
convergence to a small neighborhood of the extremum point of the static map with actuator dynamics governed by distributed diffusion PDE, which is the third contribution of this paper.

\paragraph*{Organization}
The paper is structured as follows.
The problem formulation related to the ESC with distributed diffusion PDEs
is presented in Section~\ref{sec:formulation}, including a motivating example and
the definition of the signals and system considered here.
The main results are presented in Section~\ref{sec:main-result}, 
where the distributed diffusion PDE compensation controller is designed. 
Moreover, we provide the stability analysis via Lyapunov-based  arguments and averaging theory to ensure the 
convergence of the original system into a small neighborhood of the extremum point.
Numerical results are presented in Section~\ref{sec:results}.
Finally, concluding remarks are drawn in Section~\ref{sec:conclusion}.

\paragraph*{Notation} 
The partial derivatives of a function $u(x,t)$ are denoted as 
$\partial_{x} u(x,t) = \partial u(x,t)/\partial x$, 
$\partial_{t} u(x,t) = \partial u(x,t)/\partial t$. 
Whenever convenient, we use the compact form
$u_{x}(x,t)$ e $u_{t}(x,t)$, or simply $u_{x}$ and $u_{t}$, respectively. 
The subscript ``$\mathrm{av}$'' is used to denote the average of a periodic variable with period $\Pi$. 
The Euclidean norm of a finite-dimension state variable $X(t)$ is denoted as $|X(t)|$. 
The spatial norm in ${L}_{2}[0,L]$ of the PDE state $u(x, t)$ is denoted as 
$\|{u(t)}\|^{2}_{{L}_{2}([0,L])} = \int_{0}^{L} u^2(x,t) \,dx$. The index ${L}_{2}([0,L])$ is omitted for
the sake of brevity, thus $\|{\cdot}\| = \|{\cdot}\|_{{L}_{2}([0,L])}$, if it is not specified. 
According to \cite{KH:02}, a vector-valued  function $f(t,\epsilon) \in \mathbb{R}^{n}$ 
is said of order $\mathcal{O}(\epsilon)$ over an interval $[t_{1},t_{2}]$, if 
$\exists k, \bar{\epsilon}: |{f(t,\epsilon)}| \leq k\epsilon , \forall t \in [t_{1},t_{2}]$ and 
$\forall \epsilon \in [0,\overline{\epsilon}]$. In general, the precise estimates of the constants
$k$ and $\bar{\epsilon}$ are not provided. Thus, we simply use $\mathcal{O}(\epsilon)$, 
which is interpreted as an order of magnitude relation for a sufficiently small scalar $\epsilon$.

\section{Problem Formulation}
\label{sec:formulation}

\subsection{A motivating example}

\begin{figure}[!ht]
    \centering
    \includegraphics[width=0.9\columnwidth]{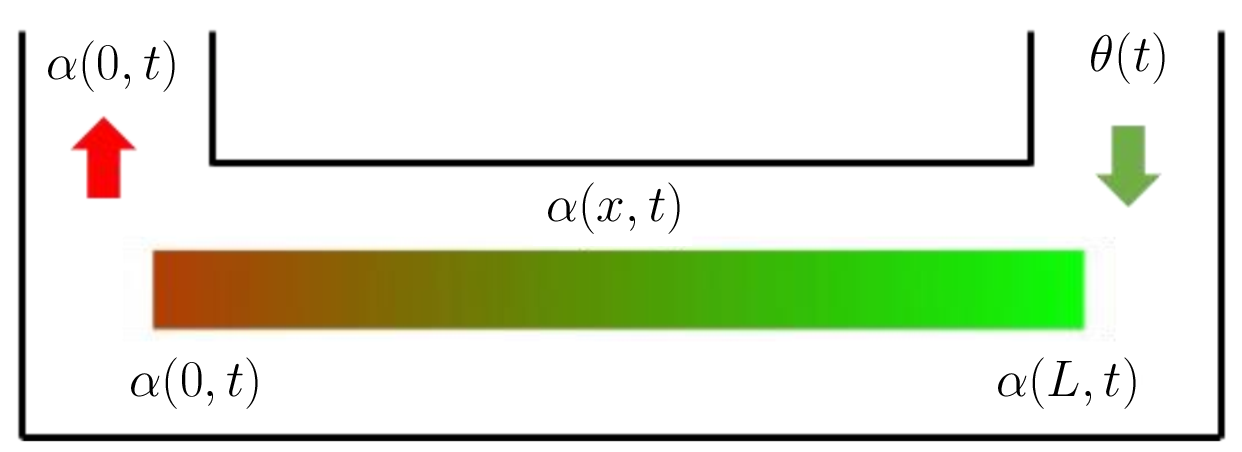}
    \caption{Geometric representation of the temperature profile modeled by the diffusion PDE in~\eqref{eq:diffusion-1}--\eqref{eq:diffusion-2}, \textcolor{black}{where $\theta(t)$ is
    the heat/cold flow source.}}
    \label{fig:temperature}
\end{figure}

Consider a diffusion (heat) PDE given by
\begin{align}
    \partial_{t}\alpha(x,t) &= \epsilon \partial_{xx}\alpha(x,t),\, \label{eq:diffusion-1}\\
    \alpha(L,t) &= \theta(t), \label{eq:diffusion-2}
\end{align}
where $\alpha(x,t)$ is the PDE state, $\epsilon > 0$ is the diffusion coefficient, 
$x \in [0,L]$, being $L$ the known domain length of the PDE.
A geometric illustration of equations~\eqref{eq:diffusion-1}--\eqref{eq:diffusion-2} 
describing the heat conduction is provided in Fig.~\ref{fig:temperature}. 

This example may represent the temperature control in a building system \cite{balaji2013zonepac}, where the heat flow from a temperature source $\theta(t)$ enters the room from the supply ventilation. 
%
Although the temperature gradient near the ceiling $\alpha(x,t)$ varies in time, it also varies along the 
space variable $x$. Thus, the heat/cold flow is diffused  from the supply to the return vent $\alpha(0,t)$.
Moreover, we can define the distributed temperature along the space variable $x$ as follows
\begin{align}
   \Theta(t) = \int_0^L \alpha(y,t) dy.
\end{align}

Thus, if we consider an optimum operating point along the space variable as $\Theta^{\ast}$,
the control problem can be formulated by designing a 
controller such that the variable $\Theta(t)$ converges towards
a small neighborhood of the desired average temperature
$\Theta^{\ast}$. This example illustrates one potential application of 
ESC for distributed diffusion PDEs, which consists of driving the variable $\Theta(t)$ towards $\Theta^{\ast}$.

\subsection{The standard gradient ESC without PDEs}
\label{ESC}

\begin{figure}[!ht]
    \centering
    \includegraphics[width=0.8\columnwidth]{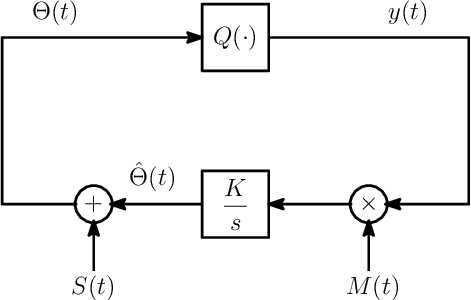}
    \caption{Standard ESC scheme without PDEs.}
    \label{fig:esc}
\end{figure}

Consider the standard gradient ESC without PDEs in Fig.~\ref{fig:esc}.
The aim is to optimize an unknown static map $y\!=\!Q(\Theta)$
by solving a real-time optimization with unknown output $y^{\ast}$ and unknown optimizer $\Theta^{*}$,
with measured output $y(t) \in \mathbb{R}$ and input~$\Theta(t) \in \mathbb{R}$.

By following the ESC methodology proposed by~\cite{Krsti2000StabilityOE},
the sinusoidal perturbation signal and the demodulation signal are given, respectively, by
\begin{equation}
S(t) = a\sin{(\omega t)}\;\; \text{and} \;\; M(t) = \frac{2}{a}\sin{(\omega t)},
\label{eq:ditherSM}
\end{equation}
where $a>0$ is the amplitude and $\omega>0$ is the angular frequency.
Both signals are judiciously designed to obtain the estimate of the unknown gradient
$\partial Q(\Theta)/\partial \Theta$ and the negative Hessian
$H:=\partial^2 Q(\Theta)/\partial \Theta^2 <0$ of the nonlinear map $Q(\Theta)$ to be maximized.
The input $\Theta(t)=\hat{\Theta}(t)+S(t)$ is derived from its real-time estimate $\hat{\Theta}(t)$ of $\Theta^{\ast}$, which is perturbed by the signal~$S(t)$. 
In this case, the estimate $\hat{\Theta}(t)$ is obtained by integrating  $\dot{\hat{\Theta}}(t) = {K}M(t)y(t)$, which locally approximates the
gradient update law $\dot{ \hat{\Theta}}(t)= KH (\hat{\Theta}(t)-\Theta^{\ast})$
by driving $\hat{\Theta}(t)$ towards $\Theta^{\ast}$.
Therefore, by defining the estimation error 
${\vartheta}(t) = \hat{\Theta}(t) - \Theta^{*}$, 
the average error dynamics can be described as 
$\dot{\vartheta}_{av } = KH{\vartheta}_{av}$, which is exponentially stable
whenever the adaptation gain is selected as $K > 0$. 

\subsection{Gradient ESC with Distributed Diffusion PDEs}

Different from the standard gradient ESC scheme shown in Fig.~\ref{fig:esc}, this paper concerns the case in which the actuation dynamics are governed by a distributed diffusion PDE.
Without loss of generality, we consider the diffusion coefficient as $\epsilon = 1$. Moreover, we consider the scalar case in which
 $\theta(t) \in \mathbb{R}$ and $\Theta \in \mathbb{R}$, such that
\begin{align}
     \Theta(t) &= \int_0^L \alpha(y,t) dy,
    \label{Theta_p_actuator}\\
        \partial_{t}\alpha(x,t) &= \partial_{xx}\alpha(x,t),\label{diffusion_actuator}\\
        \partial_{x}\alpha(0,t) &= 0,\label{boundary_actuator}\\
        \alpha(L,t) &= \theta(t),\label{theta_actuator}
\end{align}
where $\alpha:[0,L]\times\mathbb{R}_{+}\rightarrow\mathbb{R}$ and $L>0$ is the known length of the space variable.

The measured output signal is represented by the unknown static map
\begin{equation}
y(t) = Q(\Theta(t)),
\label{eq:initial_output_static_map}
\end{equation}
where $\Theta(t)$ given as in~\eqref{Theta_p_actuator} is the input of the map.

\begin{assumption}\label{aassumpstatic_map}
    The unknown static map is locally quadratic and described as follows:
    \begin{equation}
Q(\Theta) = y^{*} + \dfrac{H}{2}(\Theta - \Theta^{*})^{2},
\label{eq:static_map}
\end{equation}
\noindent
where $\Theta^{*}, \, y^{*} \in \mathbb{R}$ is the optimum output and $H<0$ is the Hessian of the map.
\end{assumption}

We emphasize that Assumption~\ref{aassumpstatic_map} is not restrictive in the sense that any $\mathcal{C}^2$ static map can be approximated as~\eqref{eq:static_map}. Thus, it follows from
\eqref{eq:initial_output_static_map} and \eqref{eq:static_map} that the output of the static map is given by
\begin{equation}
y(t) = y^{*} + \dfrac{H}{2}(\Theta(t) - \Theta^{*})^{2}.
\label{eq:final_output_static_map}
\end{equation}

\begin{figure}[!ht]
    \centering
    \includegraphics[width=\columnwidth]{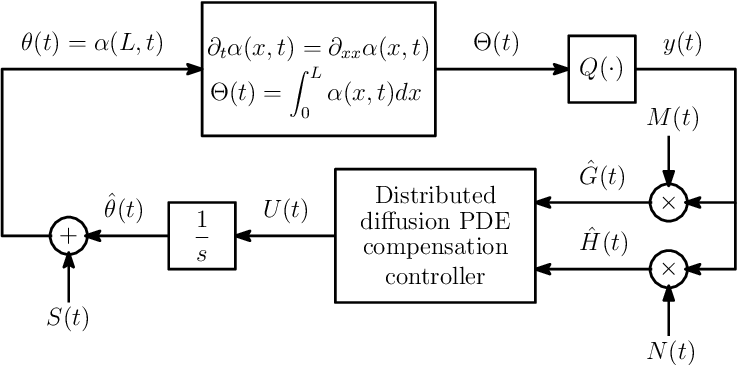}
    \caption{ESC scheme with distributed diffusion PDE governing the actuator dynamics.}
    \label{fig:esc_pde}
\end{figure}

Motivated by the ESC scheme proposed by~\cite{Oliveira2017ExtremumSF}, the ESC system with distributed diffusion PDE can be formulated considering \eqref{Theta_p_actuator}--\eqref{theta_actuator} in the ESC scheme in Section~\ref{ESC}. Therefore, the ESC scheme studied in this paper is illustrated in Fig.~\ref{fig:esc_pde}.

\subsection{Trajectory Generation to the Additive Perturbation Signal}

In the basic ESC scheme, the additive perturbation signal $S(t)$ defined as in~\eqref{eq:ditherSM} introduces the sinusoidal term $a\sin{(\omega t)}$ in the input of the static map. However, this simple choice of $S(t)$ can not introduce the appropriate sinusoidal term in the input of the static map in Fig.~\ref{fig:esc_pde} due to the presence of the distributed diffusion PDE governing the actuator dynamics.

An adequate additive perturbation signal $S(t)$ can be properly designed by solving the so-called trajectory generation problem~\cite[Chapter~12]{Krsti2008BoundaryCO}. 
Based on this reasoning, the design of the signal $S(t)$ can be performed regarding the following trajectory generation problem:
\begin{align}
    S(t) &= \beta(L,t), \label{perturbation_beta}\\
    \partial_{t}\beta(x,t) &= \partial_{xx}\beta(x,t),\label{wavekv_beta}\\
    \partial_{x}\beta(0,t) &= 0,\label{boundary_beta}\\
    \int_0^L \beta(y,t) dy &= a\sin{(\omega t)},\label{initial_cond_beta}
\end{align}
where $\beta:[0,L]\times\mathbb{R}_{+}\rightarrow\mathbb{R}$,
whose solution consists in determining
$\beta(L,t)$ such that equations~\eqref{wavekv_beta}--\eqref{initial_cond_beta}
are satisfied. 

The solution to this problem is stated in the following lemma.
\begin{lemma}\label{lem:Sdesign}
    The solution of the trajectory generation problem~\eqref{perturbation_beta}--\eqref{initial_cond_beta} is
\begin{align}
\begin{aligned}\label{eq:reference_trajectory_OMG}
    S(t) &= \frac{A}{2}e^{\sqrt{\frac{\omega}{2}}L}
    \sin{\left(\omega t + \phi + \sqrt{\frac{\omega}{2}}L\right)} \\
    &+
    \frac{A}{2}e^{-\sqrt{\frac{\omega}{2}}L}
    \sin{\left(\omega t + \phi - \sqrt{\frac{\omega}{2}}L\right)},
\end{aligned}
\end{align}
where 
\begin{gather}
    A =  2a\sqrt{\omega}/B \quad \mathrm{and} \quad
    \phi = - \psi, \label{eq:sinusoidal_conditions}
\end{gather}
and
\begin{align}
    B &= \left[e^{L\sqrt{2\omega}} + e^{-L\sqrt{2\omega}} + 2\cos{(L\sqrt{2\omega})}\right]^{1/2} \label{eq:B}\\
    \psi &= \begin{cases}
        \sign{(\psi_1)} \frac{\pi}{2}, & \mathrm{if}~\psi_2 = 0, \\
        \arctan{\left(\frac{\psi_1}{\psi_2}\right)}, & \mathrm{if}~\psi_2 > 0, \\
        \pi + \arctan{\left(\frac{\psi_1}{\psi_2}\right)}, & \mathrm{if}~\psi_2 < 0,
    \end{cases} \label{eq:psi}
\end{align}
with 
\begin{align*}
    \psi_1 &= e^{L\sqrt{\frac{\omega}{2}}}\sin{\left(L\sqrt{\tfrac{\omega}{2}}-\tfrac{\pi}{4}\right)} + 
    e^{-L\sqrt{\frac{\omega}{2}}}\sin{\left(-L\sqrt{\tfrac{\omega}{2}}-\tfrac{\pi}{4}\right)}, \\
    \psi_2 &= e^{L\sqrt{\frac{\omega}{2}}}\cos{\left(L\sqrt{\tfrac{\omega}{2}}-\tfrac{\pi}{4}\right)} + 
    e^{-L\sqrt{\frac{\omega}{2}}}\cos{\left(-L\sqrt{\tfrac{\omega}{2}}-\tfrac{\pi}{4}\right)}.
\end{align*}
\end{lemma}
\begin{IEEEproof}
Consider that the reference output has the following form
\begin{align}
\begin{aligned}\label{eq:reference_trajectory}
\beta^r(0,t) &= A \sin{(\omega t + \phi)}.
\end{aligned}
\end{align}
Since $\sin{(\omega t + \phi)} = \mathrm{Im}\left\lbrace e^{j(\omega t + \phi)}\right\rbrace$, by following the same
steps as in~\cite[Chapter~12]{Krsti2008BoundaryCO},
it is possible to show that the conditions
\eqref{perturbation_beta}--\eqref{boundary_beta}
are satisfied with
\begin{align}\label{eq:proof-lemma-1}
    \beta^r(x,t) 
    &= \frac{A}{2}e^{\sqrt{\frac{\omega}{2}}x}
    \sin{\left(\omega t + \phi + \sqrt{\frac{\omega}{2}}x\right)} \nonumber \\
    &+
    \frac{A}{2}e^{-\sqrt{\frac{\omega}{2}}x}
    \sin{\left(\omega t + \phi - \sqrt{\frac{\omega}{2}}x\right)}.
\end{align}
Thus, $S(t) = \beta^r(L,t)$ is given as in~\eqref{eq:reference_trajectory_OMG}.

Now, we show that the condition~\eqref{initial_cond_beta} is satisfied by taking $A$ and $\phi$ as in~\eqref{eq:sinusoidal_conditions}.
By integrating~\eqref{eq:proof-lemma-1} with respect to $x$ and performing some simplifications, one obtains
\begin{align}
    \int_0^L \beta^r(x,t) dx &=
    \frac{A}{2\sqrt{\omega}}
    \left[
        A_1 \sin{\left( \omega t + \phi + \phi_1 \right)} \right. \nonumber \\
        &+
        \left.
        A_2 \sin{\left( \omega t + \phi + \phi_2 \right)}
    \right],
\end{align}
where $A_1 = e^{L\sqrt{\frac{\omega}{2}}}$, $A_2 = -e^{-L\sqrt{\frac{\omega}{2}}}$, 
$\phi_1 = L\sqrt{\frac{\omega}{2}} - \frac{\pi}{4}$ and
$\phi_2 = - L\sqrt{\frac{\omega}{2}} - \frac{\pi}{4}$.
Since 
\begin{align*}
    A_1\sin{(\omega t + \phi + \phi_1)} = \mathrm{Im}\left\lbrace A_1 e^{j\phi_1}e^{j(\omega t + \phi)} \right\rbrace, \\
    A_2\sin{(\omega t + \phi + \phi_2)} = \mathrm{Im}\left\lbrace A_2 e^{j\phi_2}e^{j(\omega t + \phi)} \right\rbrace,
\end{align*}
then, 
$A_1\sin{(\omega t + \phi + \phi_1)} + 
    A_2\sin{(\omega t + \phi + \phi_2)} =
\mathrm{Im}\left\lbrace \left(A_1 e^{j\phi_1}+A_2 e^{j\phi_2}\right)e^{j(\omega t + \phi)} \right\rbrace =
\mathrm{Im}\left\lbrace B e^{j\psi}e^{j(\omega t + \phi)} \right\rbrace = B \sin{(\omega t + \phi + \psi)}$, with $B$ and $\psi$ given as in~\eqref{eq:B} and \eqref{eq:psi}, respectively. By substituting
this result in~\eqref{eq:proof-lemma-1}, it follows that
\begin{align*}
    \int_0^L \beta^r(x,t) dx = \frac{A}{2\sqrt{\omega}} B \sin{(\omega t + \phi + \psi)}.
\end{align*}
Thus, by taking $A$ and $\phi$ as in~\eqref{eq:sinusoidal_conditions}, we ensure that
\eqref{initial_cond_beta} holds. This completes the proof.
\end{IEEEproof}

Besides the additive perturbation signal $S(t)$,
we also need to specify the demodulation signals 
$M(t)$ and $N(t)$ that are used to estimate the gradient and the Hessian of the static map, respectively, by multiplying them with the output signal 
$y(t)$. Similar to \cite{Ghaffari2011MultivariableNE}, 
we consider:
\begin{equation}
\hat{G}(t) = M(t)y(t)\;\; \text{with} \;\; M(t) = \frac{2}{a}\sin{(\omega t)},
\label{eq:gradient_USOpen}
\end{equation}
and
\begin{equation}
\hat{H}(t) = N(t)y(t)\;\; \text{with} \;\; N(t) = -\dfrac{8}{a^{2}}\cos{(2\omega t)}.
\label{eq:hessian}
\end{equation}

\subsection{Estimation Error Dynamics}

Consider the following estimate variables
\begin{equation}
\hat{\theta}(t) = \theta(t) - S(t),\;\;\; \hat{\Theta}(t) = \Theta(t) - a\sin{(\omega t)},
\label{eq:estimated}
\end{equation}
and the following estimation errors
\begin{equation}
\tilde{\theta}(t) = \hat{\theta}(t) - \Theta^{*},\;\;\; \vartheta(t) =\hat{\Theta}(t) - \Theta^{*}.
\label{eq:estimated_errors}
\end{equation}

Let $\bar{\alpha}:[0,L]\times\mathbb{R}_{+}\rightarrow\mathbb{R}$ be defined as $\bar{\alpha}(x,t) = \alpha(x,t) - \beta(x,t) - \Theta^{\ast}$. 
From \eqref{Theta_p_actuator}--\eqref{theta_actuator} and \eqref{perturbation_beta}--\eqref{initial_cond_beta}, together with \eqref{eq:estimated}--\eqref{eq:estimated_errors}, it follows that:
\begin{align}
& \vartheta(t) = \int_0^L \bar{\alpha}(y,t) dy + (L-1)\Theta^{\ast}, \label{dynamics_baractuator}\\
&    \partial_{t}\bar{\alpha}(x,t) = \partial_{xx}\bar{\alpha}(x,t),\label{wavekv_baractuator}\\
&    \partial_{x}\bar{\alpha}(0,t) = 0,\label{boundary_baractuator}\\
&   \bar{\alpha}(L,t) = \tilde{\theta}(t).\label{theta_baractuator}
\end{align}

Then, the error dynamics can be obtained by taking the time derivative of \eqref{dynamics_baractuator}--\eqref{theta_baractuator}, which results in the following PDE-ODE cascade system:
\begin{align}
&    \dot{\vartheta}(t) = \int_0^L u(y,t) dy, \label{dynamics_error}\\
&    \partial_{t}u(x,t) = \partial_{xx}u(x,t),\label{wavekv_error}\\
&    \partial_{x}u(0,t) = 0,\label{boundary_actuator_error}\\
&    u(L,t) = U(t),\label{eq:30}
\end{align}
where $\dot{\hat{\theta}}=\dot{\tilde{\theta}}(t)\!\coloneqq\!U(t)$ and $u(x,t)\!\coloneqq\! \bar{\alpha}_t(x,t)$.

Note also that the propagated estimation error $\vartheta(t)$, the input $\Theta(t)$, and the optimizer of the static map $\Theta^{\ast}$ satisfy 
\begin{align}\label{eq:complete_relation}
    \vartheta(t) + a\sin{(\omega t)} = \Theta(t) - \Theta^{\ast}.
\end{align}

\section{Main Results}
\label{sec:main-result}

The main results of this paper are presented in this section. In Section~\ref{sec:C_design}, we design the controller to compensate for the distributed diffusion PDE. In Section~\ref{sec:Cav_design}, the average version of the designed compensation controller is introduced, from which a control law that can be implemented in real-time is obtained. Finally, in Section~\ref{sec:stability_proof}, the stability proof of the closed-loop system with the modified compensation controller is derived.

\subsection{Distributed Diffusion PDE Compensation Controller}
\label{sec:C_design}

The exponential stability of the 
system \eqref{dynamics_error}--\eqref{eq:30}
with the proposed controller is proved
in the following theorem.

\begin{theorem}\label{thm:stability}
Consider the closed-loop system composed by the interconnection between the PDE-ODE cascade system ~\eqref{dynamics_error}--\eqref{eq:30}
with the following controller
\begin{equation}
 U(t) = \bar{K} Z(t), \label{eq:controlador}
\end{equation}
where
\begin{gather}
  Z(t)= \vartheta(t) +\int_0^{L} g(y)u(y,t)dy,  \label{eq:transformacao-Z} \\
g(x)=\frac{1}{2}\left(L^2 - x^2 \right). \label{eq:g_definition}
\end{gather}
If $\bar{K}$ is selected such that
\begin{align}\label{eq:Kbar}
    \bar{K}< 0, \quad \mathrm{and} \quad \bar{K} \neq - (2\kappa+1)^2\pi^2/4 L^3, \quad \kappa \in \mathbb{N},
\end{align}
then, for any initial condition $u (\cdot,0)\in L^2(0,L)$, the closed-loop system has a unique exponentially stable solution given by
$(\vartheta(t),u(\cdot,t)) \in C ([0,\infty],{R}^n \times L^2 (0,L) )$.
Thus, there exist positive constants $\eta$ and $\nu$ such that
\begin{equation}\label{eq:exponential_convergence}
     \Omega(t) \leq \eta \Omega (0) e^{-\nu t},
\end{equation}
where
\begin{equation}\label{eq:augmented_state}
    \Omega(t)= |\vartheta(t)|^2 + \int_0^{L} u^2(x,t)dx.
\end{equation}
\end{theorem}
\begin{IEEEproof}
The proof follows similar steps as in~\cite{bekiaris2011compensating}.
Consider the transformation of the finite-dimensional states $\vartheta(t)$, $Z(t)$
given in~\eqref{eq:transformacao-Z} and the backstepping-like transformation \cite{bekiaris2011compensating} of the infinite-dimensional actuator state $u(x,t)$ given by:
\begin{align}
    w(x,t)&=u(x,t) - \gamma(x) \bigg(\vartheta(t) +\int_0^{L}g(y)u(y,t)dy \bigg) \nonumber \\
     &= u(x,t) - \gamma(x) Z(t), \label{eq:w_equation}
\end{align} 
where the {kernel} function $\gamma(\cdot)$ is designed as follows:
\begin{equation}
\gamma(x)= {\bar{K}} \Lambda^{-1}\left(e^{\sqrt{A_\mathrm{cl}}x} + e^{-\sqrt{A_\mathrm{cl}}x} \right),
\label{eqref:gamma_x}
\end{equation}
where
\begin{align}\label{eq:lambda}
    \Lambda =  e^{\sqrt{A_\mathrm{cl}}L} + e^{-\sqrt{A_\mathrm{cl}}L}.
\end{align}

Then, using \eqref{eq:transformacao-Z} and \eqref{eq:w_equation}, the system \eqref{dynamics_error}--\eqref{eq:30} is mapped to the following target system
\begin{align}
    \dot{Z}(t)&=A_{\mathrm{cl}} Z(t), \label{eq:Zcl} \\
   \partial_t w(x,t) &= \partial_{xx} w(x,t), \label{eq:boundary-w-1} \\
   \partial_x w(0,t) &= 0, \label{eq:boundary-w-2} \\
   w(L,t)&=0,    \label{eq:boundary-w-3}
\end{align}
where $A_{\mathrm{cl}} =\bar{K} L$.

Notice that the kernel function $\gamma(\cdot)$
given as in (\ref{eqref:gamma_x}) is well-defined for $\Lambda \neq 0$.
To show that this condition holds,
provided that $\bar{K} < 0$, it follows that $A_{\mathrm{cl}} = \bar{K} L < 0$. Let $\lambda = -\bar{K} L > 0$ such that we can write
\begin{align}\label{eq:Acl_lambda}
    A_{\mathrm{cl}} = -\lambda.
\end{align}
Thus, it follows from \eqref{eq:lambda} that
$e^{i\sqrt{\lambda}L} {+} e^{-i\sqrt{\lambda}L} = e^{i2\sqrt{\lambda}L} {+} 1$, 
which is non-zero
when $\cos{(2\sqrt{\lambda} L)} \neq -1$ and $\sin{(2\sqrt{\lambda}L)}\neq 0$, that is,
$\sqrt{\lambda} L \neq \frac{(2\kappa +1)\pi}{2}, \quad \kappa \in \mathbb{N}$.
Hence, $\lambda \neq (2\kappa+1)^2\pi^2/4L^2$, for $\kappa \in \mathbb{N}$, or $A_{\mathrm{cl}} \neq -(2\kappa+1)^2\pi^2/4L^2$, for $\kappa \in \mathbb{N}$. This is true for $\bar{K}$ selected as in~\eqref{eq:Kbar}.

The inverse transformation of \eqref{eq:boundary-w-2}--\eqref{eq:boundary-w-3} is given by
\begin{align}\label{eq:u_transformation}
    u(x,t) = w(x,t) + \gamma(x)Z(t),
\end{align}
where $ \gamma(\cdot) $ is given by~\eqref{eqref:gamma_x}.
Finally, the inverse transformation of \eqref{eq:transformacao-Z} is given by:  
\begin{align}\label{eq:inverse_transformation}
    &\vartheta(t) = \bigg(1 - \int_{0}^{L} g(y)\gamma(y) dy \bigg) Z(t)
    - \int_0^{L} g(y) w(y,t) dy.
\end{align}

Consider the following Lyapunov functional candidate
\begin{align}
   V(t) = \frac{1}{2}Z^2(t) + \frac{1}{2} \int_0^{L} w(x,t)^2dx. 
\end{align}

From~\eqref{eq:boundary-w-1}, \eqref{eq:Acl_lambda} and using integration by parts, the time derivative of the functional along the trajectories of the transformed system is
\begin{multline}
    \dot{V}(t) \leq -\lambda Z^2(t) + w(x,t)\partial_x w(x,t)|_{x=0}^{x=L} \\
    -\int_0^{L} \partial_xw(x,t)^2dx.
    \label{eqref:Vav}
\end{multline}
Using the boundary conditions \eqref{eq:boundary-w-2}, \eqref{eq:boundary-w-3}
and Poincar\'{e} inequality, we obtain
\begin{align*}
    \dot{V}(t) \leq -\lambda Z^2(t) - \frac{1}{4L^2} \int_0^{L} w(x,t)^2dx.
\end{align*}
Therefore, if we take
$\rho = \min \bigl\{-\lambda/2, 1/2L^2
\bigl\}$ we get
$\dot{V}(t) \leq - \rho V(t)$. 
Then, from the Comparison Lemma \cite{KH:02}, one has that
${V}(t) \leq V(0) e^{-\rho t}$. 

To show the stability in the original variables $\vartheta(t)$ and $u(x,t)$, it is sufficient to show that
\begin{equation}\label{eq:V_bounds}
    \underline {M}\Omega(t)\leq V(t) \leq \overline{M}\Omega(t),
\end{equation}
for positive scalars $\overline {M} $ and $ \underline {M} $ and that the condition in~\eqref{eq:exponential_convergence} is satisfied with $\eta = \overline {M}/\underline {M}$.

To obtain the bounds in~\eqref{eq:V_bounds}, 
we use
\eqref{eq:transformacao-Z}, \eqref{eq:w_equation}, \eqref{eq:u_transformation}, and \eqref{eq:inverse_transformation}, 
and we apply the Young and Cauchy-Schwartz inequalities such that the bounds in~\eqref{eq:V_bounds} are
$$ \overline {M} = \bigg( 3\omega + 9L\sup_{y\in[0,L]} \gamma(y)^2 \omega \bigg)r,$$
where $r = 1/2$, 
$\omega = 1 + L \sup_{y\in [0,L]}g(y)^2$,
$\underline{M} =1/2M$, being
$M = 9 \bigg(1 + L \sup_{y\in[0,L]}g(y)^2\gamma(y)^2\bigg)
+6L \sup_{y\in [0,L]}g_i(y)^2
+ 3 + 9 L \sup_{y\in [0,L]} \gamma (y)^2$.
Then, it follows from~\eqref{eq:Zcl} that $Z(t)$ is bounded and converges exponentially to zero. 
From~\eqref{eq:w_equation}, we can conclude that $w(\cdot,0) \in L^2(0,L)$, 
hence, it follows from~\eqref{eq:boundary-w-1}--\eqref{eq:boundary-w-3} that $w(\cdot,t) \in C (L^2(0,L))$. 
Using the inverse transformation~\eqref{eq:u_transformation},
we ensure that $u(\cdot,t) \in C (L^2(0,L))$. 

The uniqueness of the weak solution is proved using the uniqueness of the weak solution of the boundary problems
\eqref{eq:boundary-w-1}--\eqref{eq:boundary-w-3} (see~\cite{bekiaris2011compensating,evans2022partial}). This concludes the proof. 
\end{IEEEproof}

\subsection{Controller Design for the Average System}
\label{sec:Cav_design}

Although the controller in~\eqref{eq:controlador} has been effectively designed to compensate for the distributed diffusion PDE dynamics,
this controller can not be directly implemented in real-time because it depends on $\vartheta$, which is not available for measurement.
To deal with this issue, we adopt the
result in~\cite{Ghaffari2011MultivariableNE}
to the ESC system with distributed diffusion PDE. With this extension, we develop a real-time implementable version of the controller~\eqref{eq:controlador}.

By substituting~\eqref{eq:complete_relation} in~\eqref{eq:final_output_static_map}, we obtain
\begin{align}\label{eq:static_map_vartheta}
    y(t) = y^{\ast} + \frac{H}{2}(\vartheta(t)+a\sin{(\omega t)})^2.
\end{align}
Then, by substituting~\eqref{eq:static_map_vartheta} in~\eqref{eq:gradient_USOpen}, we have
\begin{align}
    \hat{G}(t) &= \frac{2}{a}\sin{(\omega t)}\left( y^{\ast} + \frac{H}{2}(\vartheta(t)+a\sin{(\omega t)})^2\right) \nonumber \\
        &= \frac{2}{a}y^{\ast}\sin{(\omega t)} + \frac{H}{a}\vartheta^2(t)\sin{(\omega t}) \nonumber \\
        &+ \frac{H}{a}\left(
        2a \vartheta(t)\sin^2{(\omega t)} + a^2 \sin^3{(\omega t)}\right).
\end{align}
Finally, by substituting~\eqref{eq:static_map_vartheta} in~\eqref{eq:hessian}, it results in
\begin{align}
    &\hat{H}(t) = -\frac{8}{a}\cos{(2\omega t)}\left(y^{\ast} + \frac{H}{2}(\vartheta(t)+a\sin{(\omega t)})^2\right) \nonumber \\
    &= -\frac{8}{a^2}y^{\ast}\cos{(2\omega t)}
    - \frac{4H}{a^2}\left[\vartheta(t)\cos{(2\omega t)}\right. \nonumber \\
    &- \left.2a\vartheta(t)\sin{(\omega t)}\cos{(2\omega t)} + a^2\sin^2{(\omega t)}\cos{(2\omega t)} \right].
\end{align}

By computing the average of $\hat{G}(t)$ and $\hat{H}(t)$, over a period $2\pi/\omega$, that is,
\begin{align*}
    \hat{G}_{\mathrm{av}}(t) = \frac{\omega}{2 \pi} \int_0^{\frac{2\pi}{\omega}} \hat{G}(t), \quad 
    \mathrm{e}
    \quad 
    \hat{H}_{\mathrm{av}}(t) = \frac{\omega}{2 \pi} \int_0^{\frac{2\pi}{\omega}} \hat{H}(t),
\end{align*}
we obtain the following average versions of the gradient and the Hessian estimates, given, respectively, by
\begin{align}\label{eq:GH_av}
    \hat{G}_{\mathrm{av}}(t)  = H \vartheta_{\mathrm{av}}(t) \quad 
    \mathrm{and} \quad
    \hat{H}_{\mathrm{av}}(t) = H,
\end{align}
where $\omega$ is selected sufficiently large such that $\vartheta(t)$ is approximately constant in the averaging analysis.

Then, if we take $\bar{K} = KH$ in~\eqref{eq:controlador}, for $K>0$ (recall that $H<0$), the average version of the controller~\eqref{eq:controlador} becomes
\begin{align}\label{eq:controlador_av}
    U_{\mathrm{av}}(t) = KH\vartheta_{\mathrm{av}}(t) + KH\int_0^{L} g(y)u_{\mathrm{av}}(y,t)dy.
\end{align}
By substituting~\eqref{eq:GH_av} in~\eqref{eq:controlador_av}, we have
\begin{align}\label{eq:controlador_medio}
    U_{\mathrm{av}}(t) = K \hat{G}_{\mathrm{av}}(t) + K \hat{H}_{\mathrm{av}}(t) \int_0^{L} g(y)u_{\mathrm{av}}(y,t)dy.
\end{align}
Furthermore, by following similar steps, we obtain the following average system of~\eqref{dynamics_error}--\eqref{eq:30}:
\begin{align}
 \dot{\vartheta}_{\mathrm{av}}(t) &= 
 \int_0^{L}u_{\mathrm{av}}(y,t)dy, \label{eq:plant-1}\\
    \partial_t u_{\mathrm{av}}(x,t) &= \partial_{xx} u_{\mathrm{av}}(x,t),  \label{eq:plant-2} \\
    \partial_x u_{\mathrm{av}}(0,t) &= 0, \label{eq:plant-3} \\
    u_{\mathrm{av}}(L,t) &= U_{\mathrm{av}}(t), \label{eq:plant-4}
\end{align}
where $\vartheta_{\mathrm{av}}(t) \in \mathbb{R}$, $U_{\mathrm{av}}(t) \in \mathbb{R}$, $x\in[0,L]$. 
By invoking Theorem~\ref{thm:stability}, 
we can conclude the exponential stability of
the average system with the controller~\eqref{eq:controlador_av}, or, equivalently, \eqref{eq:controlador_medio}.

\subsection{Stability and Convergence Analysis}
\label{sec:stability_proof}

Inspired by its average version in~\eqref{eq:controlador_medio}, and
similar to~\cite{feiling2018gradient},
we propose the following filtered average-based infinite-dimensional controller:
\begin{equation}
 U(t) = \mathcal{T}\left\lbrace  K \left[ \hat{G}(t) + \hat{H}(t)\int_0^{L} g(y) u(y,t)dy \right]\right\rbrace, \label{eq:controlador_non_average}
\end{equation}
where $g(\cdot)$ is defined in~\eqref{eq:g_definition} and the low-pass filter operator is defined as
\begin{align} \label{pedro_magic_operator}
    \mathcal{T}\{\varphi(t)\} = \mathcal{L}^{-1}\left\lbrace \frac{c}{s+c} \right\rbrace \ast \varphi(t), \quad \varphi(t):\mathbb{R}_+ \rightarrow \mathbb{R},
\end{align}
with $c > 0$ being the corner frequency, $\mathcal{L}^{-1}\{\cdot\}$ is
the inverse Laplace transform operator, and $\ast$ is the convolution operator. 

Alternatively, another important aspect concerning the implementation of the controller is the
measurement of the artificial state $u(x,t)$. By definition, we have that 
$u(x,t) = \partial_t\bar{\alpha}(x,t) = \partial_t\alpha(x,t) - \partial_t\beta(x,t)$. 
By taking into account that $\partial_t\alpha(x,t) = \partial_{xx}\alpha(x,t)$ and
$\partial_t\beta(x,t) = \partial_{xx}\beta(x,t)$, and applying integration by parts 
to the integration of~$\partial_t\beta(x,t) = \partial_{xx}\beta(x,t)$ 
under the boundary conditions~\eqref{Theta_p_actuator}--\eqref{theta_actuator} and~\eqref{perturbation_beta}--\eqref{initial_cond_beta}, together with~\eqref{eq:estimated}
the controller~\eqref{eq:controlador_non_average} can be rewritten as 
\begin{equation}
 U(t) = \mathcal{T}\left\lbrace  K \left[ \hat{G}(t) + \hat{H}(t) \left(L\hat{\theta}(t) - \Theta(t) + a \sin{(\omega t)} \right) \right]\right\rbrace. \label{eq:controlador_non_average_3}
\end{equation}
In contrast to~\eqref{eq:controlador_non_average}, the expression in
\eqref{eq:controlador_non_average_3} is more efficient to implement in real-time
because it does requires measuring only $\Theta(t)$ instead of
$\partial_t \alpha(x,t) = \partial_{xx}\alpha(x,t)$.

\begin{theorem}
\label{thm:2}
Consider the control system in Fig.~\ref{fig:esc_pde}, with the control law $U(t)$ given in (\ref{eq:controlador_non_average}) or \eqref{eq:controlador_non_average_3}. 
There exists $\omega^{\ast}>0$ such that, $\forall \omega \geq \omega^{\ast}$ sufficiently large, and $K>0$, the closed-loop system \eqref{dynamics_error}--\eqref{eq:30} has a unique locally exponentially stable periodic solution in 
%
%
$t$ with period $\Pi \coloneqq 2\pi/\omega$, denoted by  $\vartheta^{\Pi}(t)$, $u^{\Pi}(x,t)$.  
This solution satisfies the following condition:  
%
%
\begin{equation}
\begin{aligned}
& {} \left(\onenorm{\vartheta^{\Pi}(t)}^{2}+\int_0^L ({u^{\Pi}(x,t)})^{2} dx \right)^{1/2} \leq \mathcal{O}(1/\omega)\,, \quad \forall t \geq 0\,.
\label{extrachapter.eq:order_period}
\end{aligned}
\end{equation}
Moreover, 
\begin{align}
&\!\!   \limsup_{t \to \infty} |{\theta(t)\!-\!\Theta^{*}|} \!=\! \mathcal{O}\left ( \frac{A}{2}e^{\sqrt{\frac{\omega}{2}}L}\!+\!\frac{A}{2}e^{-\sqrt{\frac{\omega}{2}}L}\!+\!1/\omega\right ), \label{order_S(t)}\\
&\!\!   \limsup_{t \to \infty} \onenorm{\Theta(t)-\Theta^{*}} = \mathcal{O}\left(a+1/\omega\right), \label{Order_Theta(t)}\\
&\!\!  \limsup_{t \to \infty} \onenorm{y(t)-y^{*}} = \mathcal{O}\left(a^{2}+1/\omega^{2}\right). \label{order_y(t)}
\end{align}
\end{theorem}
\begin{IEEEproof}
Based on the exponential stability ensured in
Theorem~\ref{thm:stability} and the invertibility of the 
transformation~\eqref{eq:transformacao-Z} given in \eqref{eq:inverse_transformation}, it follows from~\eqref{eq:exponential_convergence}--\eqref{eq:augmented_state} that
\begin{align} \label{Alcaraz}
& {} \onenorm{\vartheta_{\rm{av}}(t)}^{2}\!+ 
\int_0^L {u_{\rm{av}}^2(x,t)} dx \nonumber \\ 
& {}\leq  \eta e^{-\nu t} \left( \onenorm{\vartheta_{\rm{av}}(0)}^{2}\!+\!\int_0^L {u_{\rm{av}}^2(x,0)} dx\right) \,, \quad  \forall t \!\geq\! 0\,,
\end{align}
after assuming $c \to +\infty$ in (\ref{pedro_magic_operator}) for simplicity. 
Using similar arguments as in Theorem~\ref{thm:stability}, it is possible 
to conclude from \eqref{Alcaraz} that the origin of the average closed-loop system with the distributed diffusion PDE is exponentially stable. 
In the sequel, according to the average theory in infinite dimensions~\cite{Hale1990AveragingII}, for $\omega$ sufficiently large,
the closed-loop system \eqref{dynamics_error}--\eqref{eq:30}, 
with $U(t)$ in \eqref{eq:controlador_non_average}, 
has a unique exponentially stable periodic solution around its equilibrium (origin) satisfying \eqref{extrachapter.eq:order_period}.

The asymptotic convergence to a neighborhood of the
extremum point is proved by taking the absolute value
of both sides of~\eqref{eq:complete_relation}, resulting in:
\begin{equation}
    \onenorm{\Theta(t)-\Theta^{*}} = \onenorm{\vartheta(t) + a\sin{(\omega t)}}. \label{extrachapter.vasin2}
\end{equation} 
Considering~\eqref{extrachapter.vasin2} and rewriting it 
by adding and subtracting the periodic solution 
$\vartheta^\Pi(t)$, it follows that
\begin{equation}
    \onenorm{\Theta(t)-\Theta^{*}} = \onenorm{\vartheta(t)-\vartheta^\Pi(t)+\vartheta^\Pi(t)  + a\sin{(\omega t)}}. \label{extrachapter.vasin2_FONSECA}
\end{equation} 
By invoking the average theorem \cite{Hale1990AveragingII}, 
we can conclude that $\vartheta(t)\!-\!\vartheta^\Pi(t)\!\to\!0$, exponentially. Consequently,
\begin{equation}
    \limsup_{t \to \infty} \onenorm{\Theta(t)\!-\!\Theta^{*}} = \limsup_{t \to \infty} \onenorm{\vartheta^{\Pi}(t) + a\sin{(\omega t)}} \label{extrachapter.Order_Theta(t)1}.
\end{equation}
Finally, by using the relation \eqref{extrachapter.eq:order_period}, we obtain the result in~\eqref{Order_Theta(t)}.
%


Provided that $\theta(t) - \Theta^{\ast} = \tilde{\theta}(t) + S(t)$ from~\eqref{eq:estimated} and \eqref{eq:estimated_errors}, recalling that $S(t)$ is of order $\mathcal{O}\left (  \frac{A}{2}e^{\sqrt{\frac{\omega}{2}}L}\!+\!\frac{A}{2}e^{-\sqrt{\frac{\omega}{2}}L}\right )$, as shown in  \eqref{eq:reference_trajectory_OMG}, with   $\limsup\limits_{t\rightarrow \infty} |\tilde{\theta}(t)| = \mathcal{O}\left(1/\omega\right)$, we obtain the ultimate bound given in~\eqref{order_S(t)}. 


To prove the convergence of the output $y(t)$, we employ similar arguments to obtain the bound for $\Theta(t)$ by substituting~\eqref{Order_Theta(t)} in \eqref{eq:final_output_static_map},
such that
\begin{equation}
\begin{split}
\limsup_{t \to \infty} \onenorm{y(t)-y^{*}} = \limsup_{t \to \infty} \onenorm{H\vartheta^{2}(t) + Ha^{2}\sin{(\omega t)}^{2}}. \label{extrachapter.Order_y(t)1}
\end{split}
\end{equation}
Therefore, by rewriting \eqref{extrachapter.Order_y(t)1} 
in terms of $\vartheta^{\Pi}(t)$ with the help of \eqref{extrachapter.eq:order_period}, we finally obtain \eqref{order_y(t)}.
\end{IEEEproof}

\section{Numerical Results}
\label{sec:results}

\begin{figure}[!ht]
\centering
        \includegraphics[width=\columnwidth]{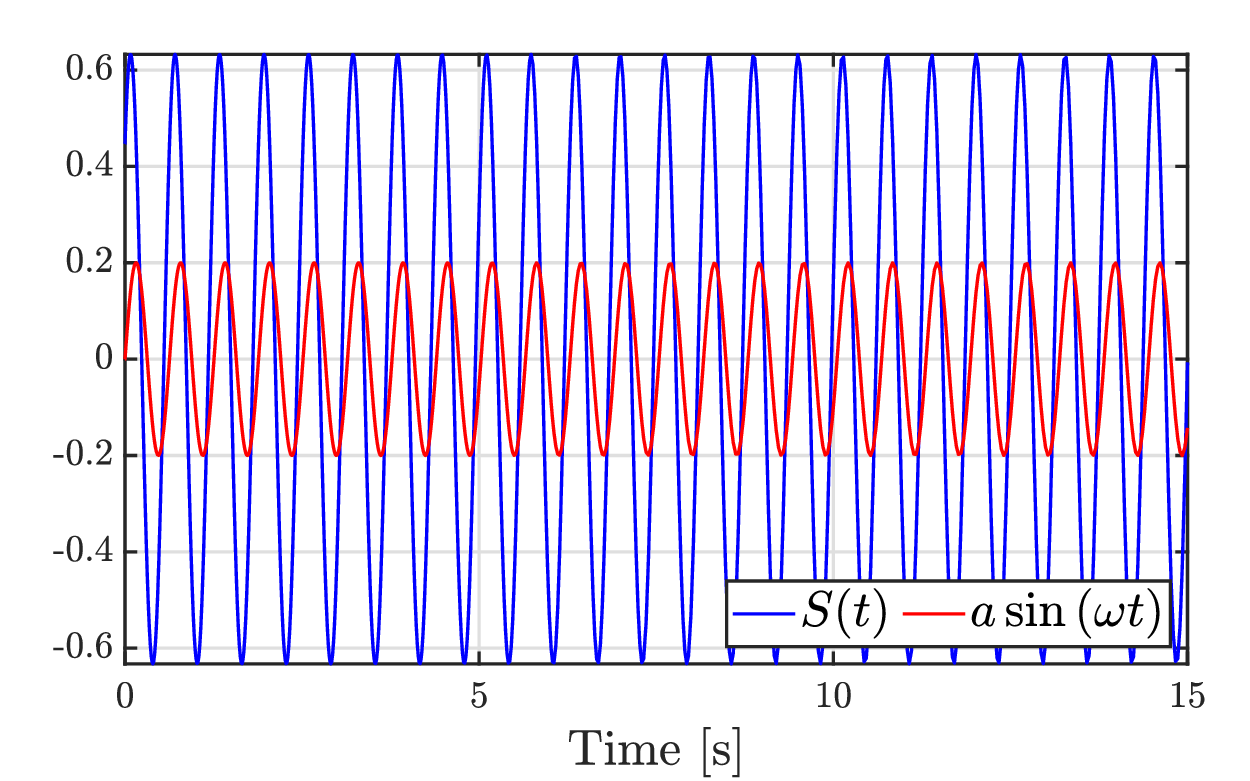}
    \caption{Illustration of the signals involved in Lemma~\ref{lem:Sdesign}.}
    \label{fig:open_loop}
\end{figure}

This section presents numerical simulations to illustrate the effectiveness of the proposed methodology.
For illustration, consider a static map~\eqref{eq:static_map}
with Hessian $H\!=\!-2$, optimizer $\Theta^{\ast}\!=\!2$ and optimal output $y^{\ast}\!=\!5$. 
For the distributed diffusion PDE in~\eqref{eq:diffusion-1}--\eqref{eq:diffusion-2} that describes
the behavior of the actuator, we consider $L = 1$ and $\epsilon=1$. 
For the controller in~\eqref{eq:controlador_non_average_3}, we consider the gain $K=0.2$ and the corner frequency 
$c = 10$ for the corresponding low-pass filter. 
For implementation purposes \cite{Ghaffari2011MultivariableNE}, we also consider a high-pass filter in the gradient estimation $\hat{G}(t)$, 
and a low-pass filter in the Hessian estimation $\hat{H}(t)$.
For the perturbation and demodulation signals, it is considered $a = 0.2$ and $\omega = 10~\mathrm{rad/s}$. 
We invoke Lemma~\ref{lem:Sdesign} to design $A$ and $\phi$ in $S(t)$ given in~\eqref{eq:reference_trajectory_OMG}. We obtain $A = 0.1356$ and $\phi = -1.4618~\mathrm{rad}$ from the provided parameters $a$ and $\omega$.
The signals $S(t)$ and $a \sin{(\omega t)}$ involved in the solution of the motion planning problem for Lemma~\ref{lem:Sdesign}
are depicted in Fig.~\ref{fig:open_loop}.

\begin{figure}[!ht]
    \centering
    \begin{subfigure}[c]{\columnwidth}
        \centering
        \includegraphics[width=\columnwidth]{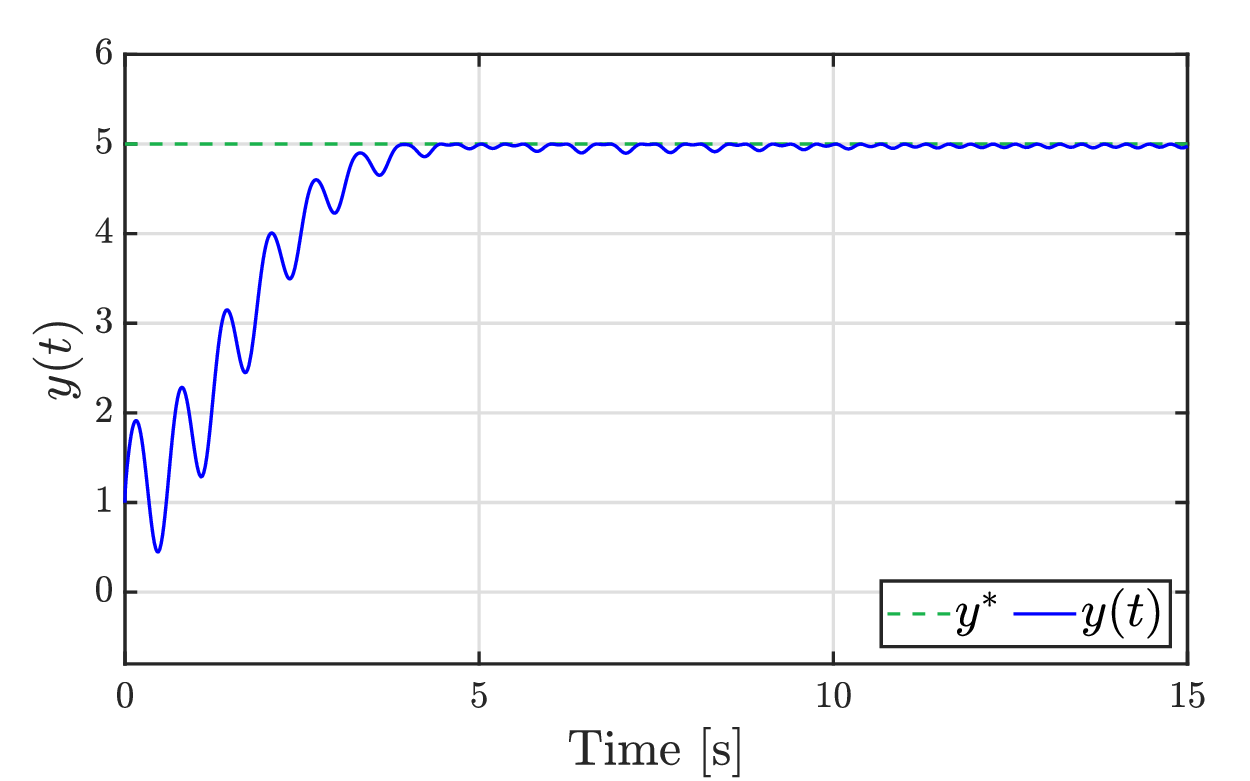}
        \caption{Output $y(t)$ of the static map.}
    \end{subfigure}
    \begin{subfigure}[c]{\columnwidth}
        \centering
        \includegraphics[width=\columnwidth]{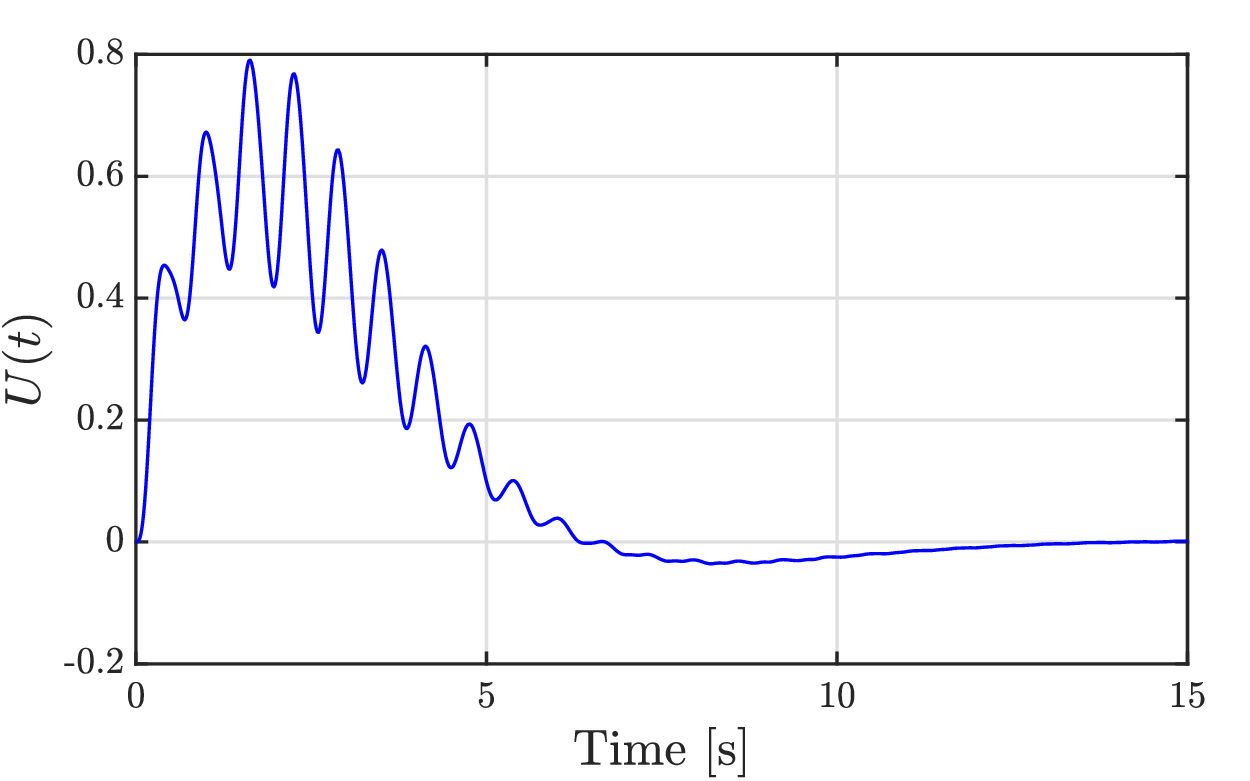}
        \caption{Control signal $U(t)$ of the controller.}
    \end{subfigure}
    \caption{Closed-loop simulation of the ESC system with actuation dynamics governed by the distributed diffusion PDE
    with the controller~\eqref{eq:controlador_non_average_3}.}
    \label{fig:closed_loop}
\end{figure}

The output of the static map (given in~\eqref{eq:static_map}) of the closed-loop system is shown in Fig.~\ref{fig:closed_loop}(a), 
where it is possible to observe the convergence to the unknown extremum point $y^{\ast} = 5$.
The signal $U(t)$ of the compensation controller implemented as in~\eqref{eq:controlador_non_average_3}  is presented in Fig.~\ref{fig:closed_loop}(b).

\begin{figure}[!ht]
    \centering
    \includegraphics[width=\columnwidth]{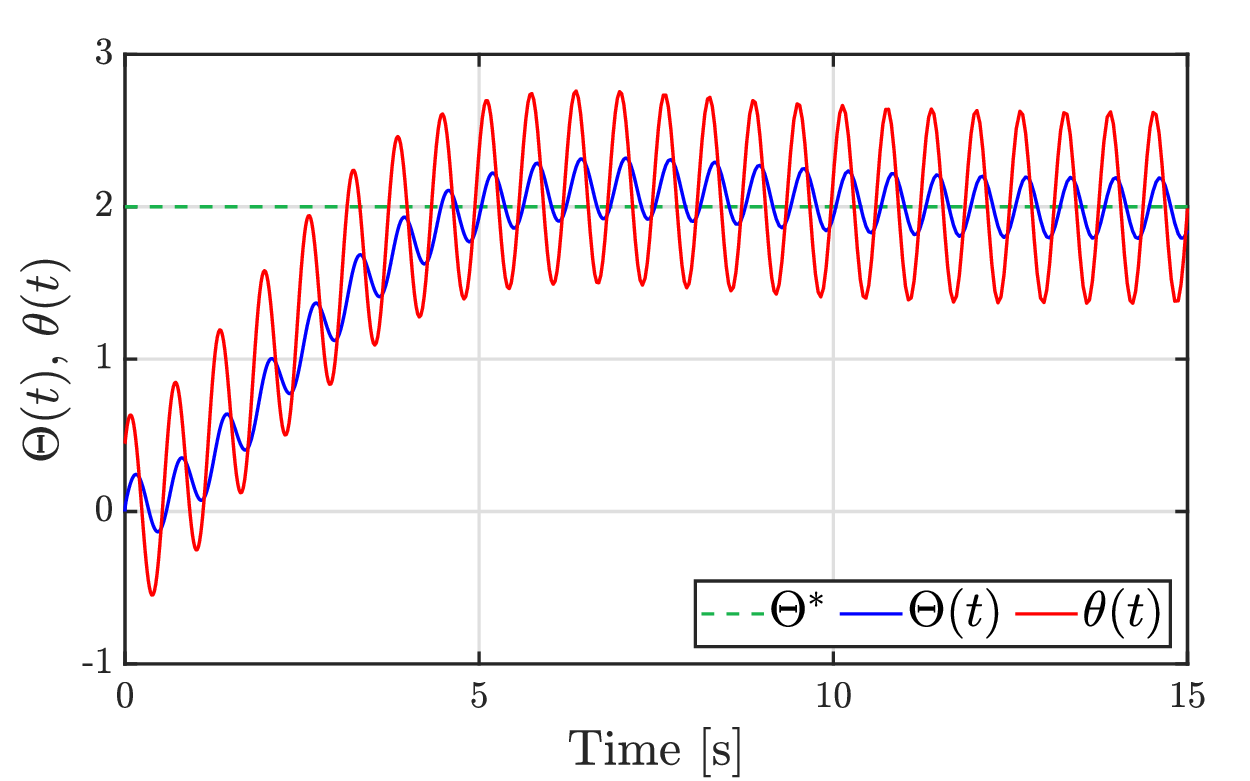}
    \caption{Signals $\theta(t)$ and $\Theta(t)$ of the closed-loop simulation.}
    \label{fig:thetas}
\end{figure}

The signals of the parameters $\theta(t)$ and $\Theta(t)$ are shown in Fig.~\ref{fig:thetas}, 
where it is also possible to note the convergence to the neighborhood of the optimal value $\Theta^{\ast} = 2$.

\begin{figure}[!ht]
    \centering
    \begin{subfigure}[b]{0.93\columnwidth}
        \centering
        \includegraphics[width=0.93\columnwidth]{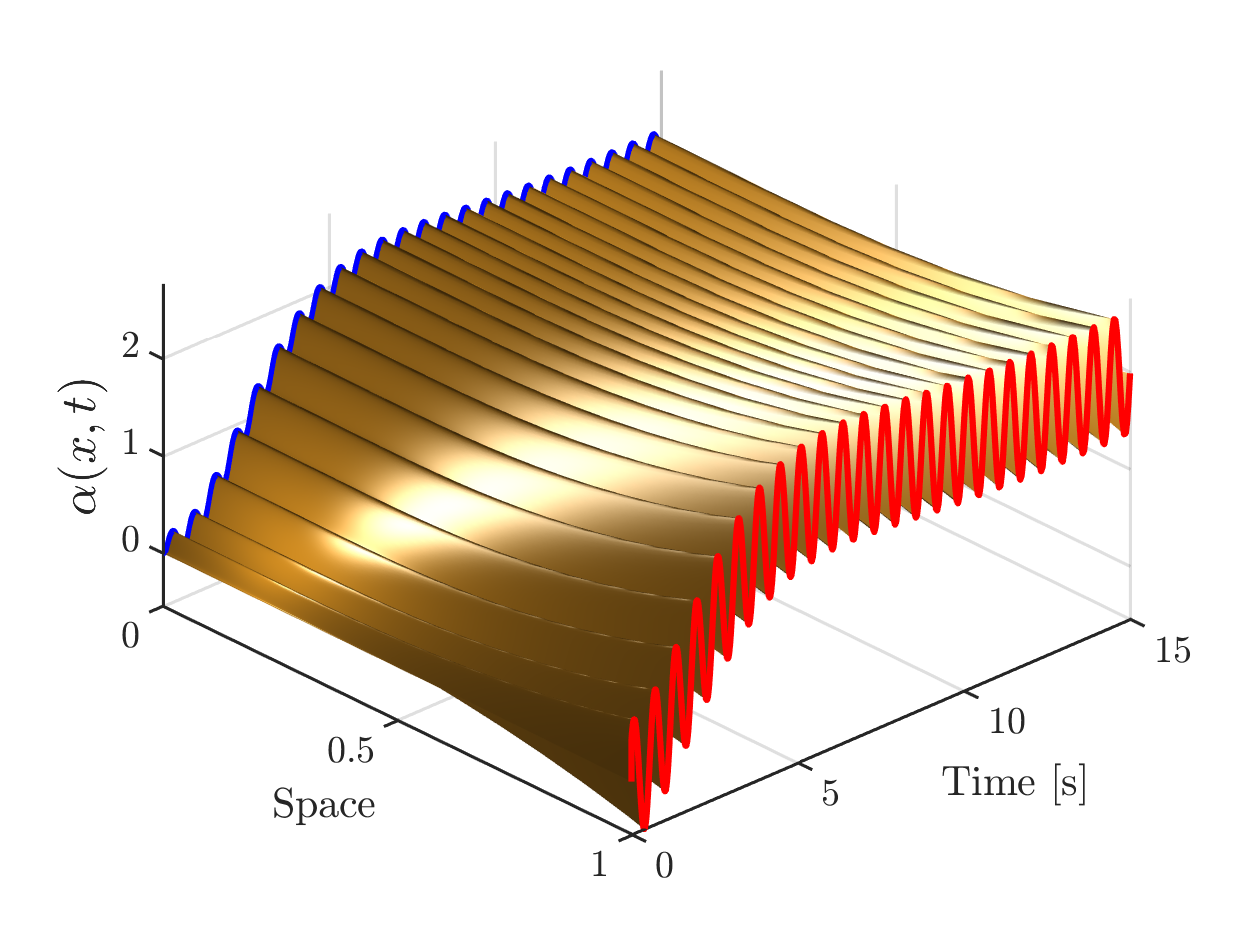}

        \caption{View of $\theta(t) = \alpha(1,t)$ in red.}
    \end{subfigure}
    \begin{subfigure}[b]{0.93\columnwidth}
        \centering
        \includegraphics[width=0.93\columnwidth]{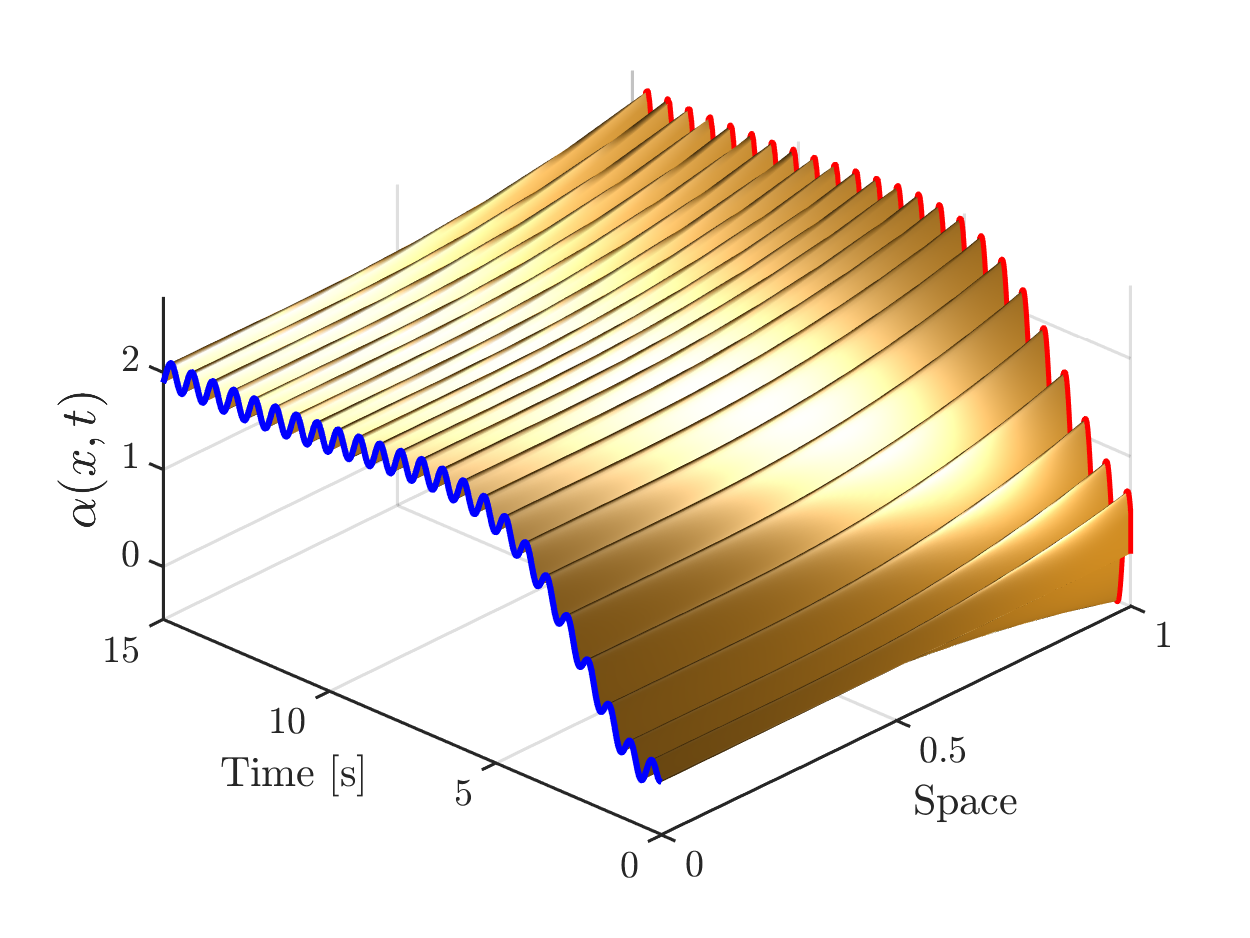}
        \caption{View of $\alpha(0,t)$ in blue.}
    \end{subfigure}
    \caption{Convergence of the state $\alpha(x,t)$ in a three-dimensional space. The signal in red is $\theta(t)=\alpha(1,t)$, while the signal in blue is $\alpha(0,t)$, both reaching a neighborhood of $\Theta^\ast = 2$.}
    \label{fig:3d-view}
\end{figure}

The evolution of the distributed diffusion PDE~\eqref{Theta_p_actuator}--\eqref{theta_actuator} of the closed-loop system in a three-dimensional space with the space domain $x \in [0,1]$ and the time $t$. The curves in blue and red show the convergence of $\alpha(0,t)$ and $\theta(t) = \alpha(1,t)$ 
to a small neighborhood around the optimizer $\Theta^{*} = 2$, respectively, as shown in Fig.~\ref{fig:3d-view}.
\vspace{-0.4cm}

\section{Conclusion}
\label{sec:conclusion}

This paper has addressed the ESC problem of locally quadratic static scalar maps with actuation dynamics governed by distributed diffusion PDEs. Concerning the appropriate design of the perturbation signal and the controller to compensate for the distributed diffusion PDE dynamics, the exponential stability of the closed-loop average system was guaranteed. Then, we have proved that the trajectories converge to a small neighborhood around the optimum point. Numerical simulations clearly illustrated the effectiveness of the proposed ESC scheme. 
\textcolor{black}{Although~\cite{Tsubakino_2023} has addressed the ESC for distributed delays, neither the control design nor the trajectory generation problem proposed in~\cite{Tsubakino_2023} could be applied to deal with distributed diffusion PDEs. Concerning the control design, Theorem~\ref{thm:stability} is different than the corresponding stability analysis in~\cite{Tsubakino_2023}, which is based on Artstein’s reduction approach. Here, we have employed a backstepping-like procedure more related to~\cite{bekiaris2011compensating}. Concerning the trajectory generation problem, the result proposed in Lemma~\ref{lem:Sdesign} is novel relative not only to~\cite{Tsubakino_2023} but also to other related work in the literature. Notably, this has been the first effort to pursue an extension of extremum seeking from the heat PDE to the class of distributed diffusion.}


\bibliographystyle{IEEEtranS}
\bibliography{IEEEabrv,references}

\end{document}